\documentclass[12pt]{article}

\usepackage{amsmath, epsfig, cite}
\usepackage{amssymb}
\usepackage{amsfonts}
\usepackage{latexsym}
\usepackage{graphicx}

\newtheorem{thm}{Theorem}[section]
\newtheorem{prop}[thm]{Proposition}

\newtheorem{exam}[thm]{\it Example}

\numberwithin{equation}{section}

\makeatletter \@addtoreset{equation}{section} \makeatother

\begin{document}
\rule{0cm}{3cm}
\begin{center}
{\Large\bf The Abel Lemma and the $q$-Gosper Algorithm}
\end{center}
 \vskip 2mm \centerline{Vincent Y. B. Chen$^1$, William Y. C. Chen$^2$, and Nancy S. S. Gu$^3$ }

\begin{center}
Center for Combinatorics, LPMC\\
Nankai University, Tianjin 300071\\
P. R. China\\

\vskip 2mm
 Email:  $^1$ybchen@mail.nankai.edu.cn, $^2$chen@nankai.edu.cn,  $^3$gu@nankai.edu.cn
\end{center}

\begin{center}
{\bf Abstract}
\end{center}

{\small Chu has recently shown that the Abel lemma on summations by
parts can serve as the underlying relation for Bailey's ${}_6\psi_6$
bilateral summation formula. In other words, the Abel lemma spells
out the telescoping nature of the ${}_6\psi_6$ sum. We present a
systematic approach to compute Abel pairs for bilateral and
unilateral basic hypergeometric summation formulas by using the
$q$-Gosper algorithm. It is demonstrated that Abel pairs can be
derived from Gosper pairs. This approach applies to many classical
summation formulas.  }

\vskip 3mm


\vskip 3mm


\section{Introduction}

We follow the notation and terminology in \cite{Gasper-Rahman}.
For $|q|<1$, the $q$-shifted factorial is defined by
$$(a;q)_\infty=
\prod_{k=0}^{\infty}(1-aq^k) \text{ and }(a;q)_n
=\frac{(a;q)_\infty}{(aq^n;q)_\infty}, \text{ for } n\in \mathbb{Z}.$$
For convenience, we shall adopt the following  notation for multiple
$q$-shifted factorials:
$$(a_1,a_2,\ldots,a_m;q)_n=(a_1;q)_n(a_2;q)_n\cdots(a_m;q)_n,$$
where $n$ is an integer or infinity.
In particular, for a nonnegative integer $k$, we have
\begin{equation} \label{x-k}
(a;q)_{-k}={1\over (aq^{-k};q)_{k}}.
\end{equation}

The (unilateral) basic hypergeometric series $_{r}\phi_s$ is
defined by
\begin{equation}
_{r}\phi_s\left[\begin{array}{cccccc} &a_1,&a_2,&\ldots,&a_r\\
&b_1,&b_2,&\ldots,&b_s
\end{array};q,z\right]=\sum_{k=0}^{\infty}\frac{(a_1,a_2,\ldots,a_r;q)_k}{(q,b_1,b_2,\ldots,b_s;q)_k}\left[(-1)^kq^{k\choose2}\right]^{1+s-r}z^k,
\end{equation}
while the bilateral basic hypergeometric series $_{r}\psi_s$ is
defined by
\begin{equation}
_{r}\psi_s\left[\begin{array}{cccccc} &a_1,&a_2,&\ldots,&a_r\\
&b_1,&b_2,&\ldots,&b_s
\end{array};q,z\right]=\sum_{k=-\infty}^{\infty}\frac{(a_1,a_2,\ldots,a_r;q)_k}{(b_1,b_2,\ldots,b_s;q)_k}\left[(-1)^kq^{k\choose2}\right]^{s-r}z^k.
\end{equation}

Recently Chu \cite{Chu} used the Abel lemma on summations by parts
to give an elementary proof of Bailey's very well-poised
$_6\psi_6$-series identity \cite{Bailey}, see also, \cite[Appendix
II.33]{Gasper-Rahman}:
\begin{eqnarray}
\label{6psi6}
\lefteqn{_6\psi_6\left[\begin{array}{cccccc} qa^{\frac{1}{2}},& -qa^{\frac{1}{2}},& b,& c,& d,& e\\[5pt]
a^{\frac{1}{2}},&- a^{\frac{1}{2}},& aq/b,& aq/c,& aq/d,& aq/e
\end{array}; q, \frac{qa^2}{bcde}\right]} \nonumber\\[10pt]
& = & \frac{(q, aq, q/a, aq/bc, aq/bd, aq/be, aq/cd, aq/ce,
aq/de;q)_{\infty}} {(aq/b, aq/c, aq/d, aq/e, q/b, q/c, q/d, q/e,
qa^2/bcde;q)_{\infty}},
\end{eqnarray}
where $|qa^2/bcde|<1$.

Let us give a brief review of Chu's approach.
The Abel lemma on summation by parts is
stated as
\begin{equation}
\sum_{k=-\infty}^{\infty}A_k(B_k-B_{k-1})=\sum_{k=-\infty}^{\infty}B_k(A_k-A_{k+1})\label{Abellemma}
\end{equation}
provided that the series on both sides are convergent.
Based on the Abel lemma, Chu found
a pair $(A_k,B_k)$:
\begin{equation}\label{6psi6Ak-c}
A_k=\frac{(b,c,d,q^2a^2/bcd;q)_k}{(aq/b,aq/c,aq/d,bcd/aq;q)_k},
\end{equation}
and
\begin{equation}\label{6psi6Bk-c}
B_k=\frac{(qe,bcd/a;q)_k}{(aq/e,q^2a^2/bcd;q)_k}\left(\frac{qa^2}{bcde}\right)^k,
\end{equation}
which leads to the following iteration relation:
\begin{eqnarray}\label{aeq}
\lefteqn{\Omega(a;b,c,d,e)=\Omega(aq;b,c,d,eq)} \nonumber \\[5pt]
&& \times
\frac{a(1-e)(1-aq)(1-aq/bc)(1-aq/bd)(1-aq/cd)}{e(1-a)(1-aq/b)(1-aq/c)(1-aq/d)(1-a^2q/bcde)},
\end{eqnarray}
where
\begin{equation}\label{6psi6l}
\Omega(a;b,c,d,e)={_6}\psi_6\left[\begin{array}{cccccc} qa^{\frac{1}{2}},& -qa^{\frac{1}{2}},& b,& c,& d,& e\\[5pt]
a^{\frac{1}{2}},& -a^{\frac{1}{2}},& aq/b,& aq/c,& aq/d,& aq/e
\end{array}; q, \frac{qa^2}{bcde} \right].
\end{equation}

Because of the symmetries in $b, c, d, e$, applying the identity
(\ref{aeq}) three times with respect to the parameters
 $a$ and $d$, $a$ and $c$, $a$ and $b$, we arrive at the following iteration
 relation:
\begin{eqnarray}
\Omega(a;b,c,d,e) &=& \Omega(aq^4;bq,cq,dq,eq)\nonumber \\[5pt]
&& \times \frac{a^4 q^6}{bcde}\cdot \frac{1-aq^4}{1-a}\cdot
\frac{(1-b)(1-c)(1-d)(1-e)}{(a^2q/bcde;q)_4} \nonumber \\[5pt]
&&\times
\frac{(aq/bc,aq/bd,aq/be,aq/cd,aq/ce,aq/de;q)_2}{(aq/b,aq/c,aq/d,aq/e;q)_3}.
\end{eqnarray}
Again, iterating the above relation $m$ times, we get
\begin{eqnarray}
\Omega(a;b,c,d,e) &=& \Omega(aq^{4m};bq^m,cq^m,dq^m,eq^m)\nonumber \\[5pt]
&& \times \frac{a^{4m}q^{6m^2}}{(bcde)^m}\cdot
\frac{1-aq^{4m}}{1-a}\cdot
\frac{(b,c,d,e;q)_m}{(a^2q/bcde;q)_{4m}} \nonumber \\[5pt]
&&\times
\frac{(aq/bc,aq/bd,aq/be,aq/cd,aq/ce,aq/de;q)_{2m}}{(aq/b,aq/c,aq/d,aq/e;q)_{3m}}.
\end{eqnarray}
Replacing the summation index $k$ with  $k-2m$, we obtain the
transformation formula
\begin{eqnarray}\label{bcde}
\lefteqn{\Omega(a;b,c,d,e) = \Omega(a;bq^{-m},cq^{-m},dq^{-m},eq^{-m})} \nonumber \\[5pt]
&&\times
\frac{(aq/bc,aq/bd,aq/be,aq/cd,aq/ce,aq/de;q)_{2m}}{(q/b,q/c,q/d,q/e,aq/b,aq/c,aq/d,aq/e;q)_{m}(a^2q/bcde;q)_{4m}}.
\end{eqnarray}
Setting $m \rightarrow \infty$, Chu obtained the $_{6}\psi_6$
summation formula \eqref{6psi6} by Jacobi's triple product identity
\cite[Appendix II.28]{Gasper-Rahman}
\begin{equation}
\sum_{k=-\infty}^{\infty}q^{k^2}z^k=(q^2,-qz,-q/z;q^2)_\infty,
\label{jacobi}
\end{equation}
since
\begin{eqnarray}
\lefteqn{\lim_{m \rightarrow \infty} \Omega(a;bq^{-m},cq^{-m},dq^{-m},eq^{-m})}\nonumber \\[5pt]
&=& \lim_{m \rightarrow \infty}
{_6\psi_6} \left[\begin{array}{cccccc} qa^{\frac{1}{2}},& -qa^{\frac{1}{2}},& bq^{-m},& cq^{-m},& dq^{-m},& eq^{-m}\\[5pt]
a^{\frac{1}{2}},& -a^{\frac{1}{2}},& aq^{m+1}/b,& aq^{m+1}/c,&
aq^{m+1}/d,& aq^{m+1}/e
\end{array}; q, \frac{a^2q^{4m+1}}{bcde} \right] \nonumber \\[5pt]
&=& \sum_{k=-\infty}^{\infty}
\frac{1-aq^{2k}}{1-a}q^{2k^2-k}a^{2k}=\frac{1}{1-a}
\sum_{k=-\infty}^{\infty}(-1)^k q^{{k \choose 2}}a^k
= (aq,q/a,q;q)_{\infty}.
\end{eqnarray}

This paper is motivated by the question of how  to systematically
compute Abel pairs for bilateral summation formulas. We find that
the $q$-Gosper is an efficient mechanism for this purpose. The
$q$-Gosper algorithm has been extensively studied.  Koorwinder
 gave a rigorous description of the $q$-Gosper
algorithm in \cite{Koornwinder}. Abramov-Paule-Petkov\v{s}ek
\cite{qgosperAP} developed the algorithm \verb"qHyper" for finding
all $q$-hypergeometric solutions of linear homogeneous recurrences
with polynomial coefficients. Later B\"{o}ing-Koepf \cite{qgosperBK}
gave an algorithm for the same purpose. The Maple package
\verb"qsum.mpl" was described by B\"{o}ing-Koepf \cite{qgosperBK}.
In \cite{Riese}, Riese presented an
 generalization of the $q$-Gosper's algorithm to
indefinite bibasic hypergeometric summations.

 Recall that a function $t_k$ is called a
basic hypergeometric term if $t_{k+1}/t_k$ is a rational function
of $q^k$. The $q$-Gosper algorithm is devised to answer  the
question if there is a basic hypergeometric term $z_k$ for a given
basic hypergeometric term $t_k$ such that
\begin{equation}t_k=z_{k+1}-z_k.\label{intro1}
\end{equation}

We observe that for an iteration relation of a summation formula,
the difference of the $k$th term of the both sides is a basic
hypergeometric term for which the $q$-Gosper algorithm can be
employed.

The main result of this paper is to present a general framework to
deal with basic hypergeometric identities based on the $q$-Gosper
algorithm. We always start with an iteration relation. Then we use
the $q$-Gosper algorithm to generate a Gosper pair $(g_k,h_k)$ if it
exists. We next turn to the iteration relation and derive the
desired identity by computing the limit value. Actually, once
aGosper pair $(g_k,h_k)$ is obtained,  one can easily compute the
corresponding Abel pair. Indeed, the Abel pair for the ${}_6\psi_6$
sum discovered by Chu \cite{Chu} agrees with the Abel pair derived
from the Gopser pair by using our approach. In general, our method
is efficient for many classical summation formulas with parameters.

As examples, we give Gosper pairs and Abel pairs of several
well-known bilateral summation formulas including Ramanujan's
${}_1\psi_1$ summation formula.  In the last section we demonstrate
that the idea of Gosper pairs can be applied to unilateral summation
formulas as well. We use the $q$-Gauss ${}_2\phi_1$ summation
formula as an example to illustrate the procedure to compute Gosper
pairs. As another example, we derive the Gosper pair and the Abel
pair of the ${}_6\phi_5$ summation formula.

In comparison with a recent approach presented by Chen-Hou-Mu
\cite{Chen-Hou-Mu} for proving nonterminating basic hypergeometric
series identities by using the $q$-Zeilberger algorithm, one sees
that the approach we undertake in this paper does not rely on the
introduction of the parameter $n$ in order to establish recurrence
relations, and only makes use of the $q$-Gosper algorithm.

\section{The Gosper Pairs for Bilateral Summations}

In spite of its innocent looking,  the Abel lemma is intrinsic
 for some sophisticated bilateral basic hypergeometric
identities. In this section, we introduce the notion of Gosper pairs
and show that one may apply the $q$-Gosper algorithm to construct
Gosper pairs which can be regarded as certificates like the Abel
pairs to justify iteration relations for bilateral summations.
Furthermore, it is easily seen that one can compute the Abel pairs
from Gosper pairs.

Suppose that we have a bilateral series $\sum_{-\infty}^{\infty}
F_k(a_1, a_2, \ldots, a_n)$ which has a closed form. Making the
substitutions $a_i\rightarrow a_iq$ or $a_i\rightarrow a_i/q$ for
some parameters $a_i$, the closed product formula induces an
iteration relation for the summation which can be stated as an
identity of the form:
\begin{equation}
\sum_{k=-\infty}^{\infty}F_k(a_1,a_2,\ldots,a_n)
=\sum_{k=-\infty}^{\infty}G_k(a_1,a_2,\ldots,a_n).\label{abel3}
\end{equation}

We assume that $\sum_{k=-\infty}^{\infty}F_k(a_1,a_2,\ldots,a_n)$
and $\sum_{k=-\infty}^{\infty}G_k(a_1,a_2,\ldots,a_n)$ are both
convergent. We also assume that
\begin{equation}\displaystyle{\lim_{k\rightarrow
\infty}}h_k=\displaystyle{\lim_{k\rightarrow -\infty}}h_k
   .
   \label{limitc}
\end{equation}
 We note that there are many bilateral summations with the above limit property.

An Gosper pair $(g_k,h_k)$ is a pair of basic hypergeometric terms
such that
\begin{eqnarray*}\label{abel4}
         g_k-h_k  & = &  F_k(a_1, a_2, \ldots, a_n), \\[5pt]
         g_k - h_{k+1} &   = &  G_k (a_1, a_2, \ldots, a_n).
\end{eqnarray*}
Evidently, once a Gosper pair is derived, the identity (\ref{abel3}) immediately
justified by the following telescoping relation:
\begin{equation}
\sum_{k=-\infty}^{\infty}(g_k-h_k)=
\sum_{k=-\infty}^{\infty}(g_k-h_{k+1}).\label{rAbel}
\end{equation}

We are now ready to describe our approach. Let us take
 Ramanujan's ${}_1\psi_1$ sum
\cite[Appendix II.29]{Gasper-Rahman} as an example:
\begin{equation}
\label{1psi1}
_1\psi_1\left[\begin{array}{c} a\\[5pt]
b
\end{array}; q, z \right] = \frac{(q, b/a, az, q/az;q)_{\infty}}
{(b,q/a,z,b/az;q)_{\infty}},
\end{equation}
where $|b/a|<|z|<1$.
There are many proofs of this identity, see, for example, Hahn \cite{Hahn}, Jackson
\cite{Jackson}, Andrews \cite{Andrews1,Andrews2}, Ismail
\cite{Ismail}, Andrews and Askey \cite{Andrews-Askey}, and Berndt
\cite{Berndt}.

\begin{prop} \label{prop1psi1}The following is a Gosper pair for
Ramanujan's ${}_1\psi_1$ sum:
\begin{eqnarray*}
g_k=\frac{az^{k+1}}{(az-b)} \frac{(a;q)_k}{(b;q)_k}, \\[5pt]
h_k=\frac{bz^k}{(az-b)} \frac{(a;q)_k}{(b;q)_k}.
\end{eqnarray*}
\end{prop}

Step 1. Construct an iteration relation from the closed product
form, namely, the right hand side of (\ref{1psi1}).

Setting $b$ to $bq$ in (\ref{1psi1}), we get
\begin{equation}
\label{1psi1bq}
_1\psi_1\left[\begin{array}{c} a\\[5pt]
bq
\end{array}; q, z \right] = \frac{(q, bq/a, az, q/az;q)_{\infty}}
{(bq,q/a,z,bq/az;q)_{\infty}}.
\end{equation}
Define
\begin{equation}
f(a,b,z)=\ _1\psi_1\left[\begin{array}{c} a\\[5pt]
b
\end{array}; q, z \right].
\end{equation}
 Comparing the right hand sides of \eqref{1psi1} and
\eqref{1psi1bq} gives the following iteration relation
(see also \cite{Andrews-Askey})
\begin{equation}\label{1psi1iterate}
_1\psi_1\left[\begin{array}{c} a\\[5pt]
b
\end{array}; q, z \right]=\frac{(1-b/a)}{(1-b)(1-b/az)}\ {_1}\psi_1\left[\begin{array}{c} a\\[5pt]
bq
\end{array}; q, z \right].
\end{equation}
Notice that both sides of the above identity are convergent.

Let $F_k(a,b,z)$ and $G_k(a,b,z)$ denote the $k$th terms of the left hand side and
the right hand side summations in \eqref{1psi1iterate} respectively, that is,
\begin{equation}
F_k(a,b,z)=\frac{(a;q)_k}{(b;q)_k}z^k \ \text{ and }\
G_k(a,b,z)=\frac{(1-b/a)}{(1-b)(1-b/az)}F_k(a,bq,z).
\end{equation}

Step 2. Apply the $q$-Gosper algorithm to try to find a Gosper pair
$(g_k,h_k)$.

It is essential to observe that $F_k(a,b,z)-G_k(a,b,z)$ is a basic
hypergeometric term. In fact it can be written as
\begin{equation*}
\left(1-bq^k-\frac{1-b/a}{1-b/az}\right)\frac{(a;q)_k}{(b;q)_{k+1}}z^k.
\end{equation*}

Now we may employ the $q$-Gosper algorithm for the following equation
\begin{equation} \label{fgkhk}
F_k(a,b,z)-G_k(a,b,z)=h_{k+1}-h_k,
\end{equation}
 and we find a solution of simple form
\begin{equation}\label{h}
h_k=\frac{bz^k}{(az-b)} \frac{(a;q)_k}{(b;q)_k},
\end{equation}
which also satisfied  the limit condition
\[ \lim_{k \rightarrow \infty} h_k = \lim_{k\rightarrow - \infty}  h_k  = 0.\]
As far as the verification of (\ref{1psi1iterate}) is concerned, the
existence of a solution $h_k$ and the limit condition (\ref{limitc})
would guarantee that the identity holds. Now it takes one more step
to compute the Gosper pair:
\begin{equation}\label{g}
g_k=h_k+F_k(a,b,z)=\frac{az^{k+1}}{(az-b)}
\frac{(a;q)_k}{(b;q)_k}.
\end{equation}

Step 3. Based on the iteration relation and the limit values,
     one can verify the summation formula.

From the iteration relation \eqref{1psi1iterate}, we may reduce
the evaluation of the bilateral series ${}_1\psi_1$ to a special
case
\begin{equation}\label{i1psi1}
f(a,b,z) = \frac{(b/a;q)_\infty} {(b,b/za;q)_\infty}\ f(a,0,z).
\end{equation}
Setting $b=q$ in \eqref{i1psi1}, we get
\begin{equation}
f(a,0,z) = \frac{(q,q/za;q)_\infty}{(q/a;q)_\infty}
\sum_{k=-\infty}^{\infty} \frac{(a;q)_k}{(q;q)_k}z^k \nonumber.
\end{equation}
Invoking the relation (\ref{x-k}), we see that $1/(q, q)_{-k}=0$ for
any positive integer $k$. Consequently, the above bilateral sum
reduces to a unilateral sum. Exploiting the $q$-binomial theorem
\cite[Appendix II.3]{Gasper-Rahman}
\begin{equation}
\sum_{k=0}^{\infty}\frac{(a;q)_k}{(q;q)_k}z^k
=\frac{(az;q)_\infty}{(z;q)_\infty},    \label{binomial}
\end{equation}
we get the evaluation
\begin{equation}
f(a,0,z) =
\frac{(q,az,q/az;q)_\infty}{(q/a,z;q)_\infty}.\label{1psi12}
\end{equation}
Hence the identity \eqref{1psi1} follows from (\ref{i1psi1}) and
(\ref{1psi12}).

It should be warned that it is not always the case that there is a
solution $h_k$ to the equation \eqref{fgkhk} in general case.
 If one encounters this
scenario, one should still have alternatives to try another
iteration relations, as is done for the ${}_3\psi_3$ sum
 in Example \ref{example}.

Let us now examine how to generate an Abel pair $(A_k, B_k)$ from a
Gosper pair $(g_k, h_k)$. Setting
\begin{equation}\label{ghab}
g_k=A_kB_k \ \ \text{and}\ \ h_k=A_kB_{k-1}, \end{equation}
then we
see that
\begin{equation}\label{ib}
\frac{B_k}{B_{k-1}}=\frac{g_k}{h_k}.
\end{equation}
Without loss of generality, we may assume that $B_0=1$.
Iterating (\ref{ib}) yields an Abel pair $(A_k,B_k)$.

For the Ramanujan's ${}_1\psi_1$ sum \eqref{1psi1},  we can compute
the Abel pair by using the $q$-Gosper algorithm.

\begin{prop} \label{aprop1psi1}
The following  is an Abel pair for
Ramanujan's ${}_1\psi_1$ sum:
\begin{eqnarray*}
A_k & = & \frac{az}{az-b}
\frac{(a;q)_k}{(b;q)_k}\left(\frac{b}{a}\right)^k
,\\[8pt]
B_k  & = & \left(  {az\over b} \right)^k .
\end{eqnarray*}
\end{prop}

It is a routine to verify $(A_k, B_k)$ is indeed an Abel pair
for the ${}_1\psi_1$ sum. First we have
\begin{eqnarray}
B_k-B_{k-1} & = & \left(\frac{az}{b}-1\right)\left(\frac{az}{b}\right)^{k-1},\\[8pt]
A_k -A_{k+1} & = & \frac{(1-b/a)}{(1-b)(1-b/az)}
  \frac{(a;q)_k}{(bq;q)_k} \left(\frac{b}{a}\right)^k.
\end{eqnarray}
Then the iteration relation \eqref{1psi1iterate}
is deduced from the Abel lemma:
\begin{eqnarray}
\lefteqn{\sum_{k=-\infty}^{\infty} \frac{az}{az-b}
\frac{(a;q)_k}{(b;q)_k} \left(\frac{b}{a}\right)^k
 \left(\frac{az}{b}-1\right)\left(\frac{az}{b}\right)^{k-1}}\\[8pt]
&=&\sum_{k=-\infty}^{\infty}\left( {az\over b} \right)^k
 \frac{(1-b/a)}{(1-b)(1-b/az)}
\frac{(a;q)_k}{(bq;q)_k} \left({b\over a} \right)^k.
\end{eqnarray}

We next give some  examples for bilateral summations.

\begin{exam}
The sum of a well-poised $_2\psi_2$ series
{\rm($\!$\cite{Gasper-Rahman}, Appendix II.30):}
\begin{equation}
\label{2psi2}
_2\psi_2\left[\begin{array}{cc} b,& c\\[5pt]
aq/b,& aq/c
\end{array}; q, -\frac{aq}{bc} \right] = \frac{(aq/bc;q)_{\infty}(aq^2/b^2,aq^2/c^2,q^2,aq,q/a;q^2)_{\infty}}
{(aq/b,aq/c,q/b,q/c,-aq/bc;q)_{\infty}} ,
\end{equation}
where $|aq/bc|<1$.
\end{exam}

Write the $k$th term of the left hand side of (\ref{2psi2}) as
\begin{equation}
F_k(a,b,c)=\frac{(b,c;q)_k}{(aq/b,aq/c;q)_k}\left(-\frac{aq}{bc}\right)^k.
\end{equation}
Substituting $b$ with $b/q$ in (\ref{2psi2}), we are led to the
iteration relation
\begin{align}\label{2psi2iterate}
\sum_{k=-\infty}^{\infty}F_k(a,b,c)&=\frac{(1-aq/bc)(1-aq^2/b^2)}{(1+aq/bc)(1-q/b)(1-aq/b)}\sum_{k=-\infty}^{\infty}F_k(a,b/q,c).
\end{align}
Let
\begin{equation}
G_k(a,b,c)=\frac{(1-aq/bc)(1-aq^2/b^2)}{(1+aq/bc)(1-q/b)(1-aq/b)}F_k(a,b/q,c).
\end{equation}
Implementing  the $q$-Gosper algorithm, we obtain a Gosper pair
\begin{eqnarray}
g_k & = & \frac{(b^2cq^k-aq^2)}{(aq+bc)(bq^k-q)}F_k(a,b,c),\\[6pt]
h_k & = & \frac{bq(c-aq^k)}{(aq+bc)(bq^k-q)}F_k(a,b,c).
\end{eqnarray}
The companion Abel pair is given below:
\begin{eqnarray}
A_k  & =  &  \frac{(b,c;q)_k(b^2cq^k-aq^2)}{(aq/b,b^2c/aq;q)_k(aq+bc)(-q+bq^k)},
\\[5pt]
B_k & = & \frac{(b^2c/aq;q)_k}{(aq/c;q)_k}\left(-\frac{aq}{bc}\right)^k.\label{2psi2akbk}
\end{eqnarray}

Noticing that \eqref{2psi2} is symmetric in $b$ and $c$, we have
\begin{align}\label{2psi2iterate2}
\sum_{k=-\infty}^{\infty}F_k(a,b,c)&=\frac{(1-aq/bc)(1-aq^2/bc)(1-aq^2/b^2)(1-aq^2/c^2)}{(1+aq/bc)(1+aq^2/bc)(1-q/b)(1-q/c)(1-aq/b)(1-aq/c)}\nonumber\\
&\quad\times\sum_{k=-\infty}^{\infty}F_k(a,b/q,c/q).
\end{align}

Finally, we can reach
 \eqref{2psi2} by iterating \eqref{2psi2iterate2} infinitely many times
 along with Jacobi's
triple product identity \eqref{jacobi} as the limit case.

\begin{exam}
Bailey's sum of a well-poised $_3\psi_3$
{\rm($\!$\cite{Gasper-Rahman}, Appendix II.31):}\label{example}
\begin{equation}
\label{3psi3}
_3\psi_3\left[\begin{array}{ccc} b,&c,&d\\[5pt]
q/b,&q/c,&q/d
\end{array}; q, \frac{q}{bcd} \right] =
\frac{(q,q/bc,q/bd,q/cd;q)_{\infty}}
{(q/b,q/c,q/d,q/bcd;q)_{\infty}},
\end{equation}
where $|q/bcd|<1.$
\end{exam}

Substituting $d$ with $d/q$ in (\ref{3psi3}), one obtains the
iteration relation
\begin{eqnarray}
\label{3psi3iterate}
&&\quad_3\psi_3\left[\begin{array}{ccc} b,&c,&d\\[5pt]
q/b,&q/c,&q/d
\end{array}; q, \frac{q}{bcd} \right]\nonumber\\
&&=\frac{(1-q/bd)(1-q/cd)}{(1-q/d)(1-q/bcd)}\ _3\psi_3\left[\begin{array}{ccc} q,&c,&d/q\\[5pt]
q/b,&q/c,&q^2/d
\end{array}; q, \frac{q^2}{bcd} \right].
\end{eqnarray}

We remark that this sum is in fact an example for which the
$q$-Gosper algorithm does not succeed for the iteration relation
derived from a straightforward
 substitution such as
 $d\rightarrow dq$ or $d \rightarrow d/q$. Instead, using an idea of Paule \cite{paule}
 of symmetrizing a bilateral summation, we replace $k$ by $-k$ on the left hand side of \eqref{3psi3} to get
\begin{equation}
\label{3psi31}
_3\psi_3\left[\begin{array}{ccc} b,&c,&d\\[5pt]
q/b,&q/c,&q/d
\end{array}; q, \frac{q^2}{bcd} \right].
\end{equation}
Let $F_k(b,c,d)$ be the average of the $k$th summands
of \eqref{3psi3} and \eqref{3psi31}, namely,
\begin{equation}
F_k(b,c,d)=
\frac{(b,c,d;q)_k}{(q/b,q/c,q/d;q)_k}\left(\frac{q}{bcd}\right)^k\frac{1+q^k}{2},
\end{equation}
and  let
\begin{equation}
G_k(b,c,d)=\frac{(1-q/bd)(1-q/cd)}{(1-q/d)(1-q/bcd)}F_k(b,c,d/q).
\end{equation}
With regard to $F_k(b,c,d)-G_k(b,c,d)$, the $q$-Gosper algorithm
generates a Gosper pair:
\begin{eqnarray}g_k & = &
\frac{bdq^{k+1}+cdq^{k+1}-bcd^2q^k-q^2+dq^{k+1}
+bcdq^{k+1}-bcd^2q^{2k}-q^{k+2}}{(1+q^k)(bcd-q)(q-dq^k)}\nonumber
\\
& & \times F_k(b,c,d)\label{3psi3xk},
\\[6pt]
h_k & = &
\frac{d(b-q^k)(c-q^k)}{(1+q^k)(q-bcd)(1-dq^{k-1})}F_k(b,c,d)\label{3psi3yk},
\end{eqnarray}
which implies the iteration relation \eqref{3psi3iterate}. Invoking
the symmetric property of the parameters $b,c$ and $d$, we have
\begin{align}
\label{3psi3iterate2}
&_3\psi_3\left[\begin{array}{ccc} b,&c,&d\\[5pt]
q/b,&q/c,&q/d
\end{array}; q, \frac{q}{bcd}
\right]=\frac{(1-q/bc)(1-q^2/bc)(1-q/bd)(1-q^2/bd)}{(1-q/b)(1-q/c)(1-q/d)(1-q/bcd)}\nonumber\\
&\quad\quad\quad\quad\quad\times\frac{(1-q/cd)(1-q^2/cd)}{(1-q^2/bcd)(1-q^3/bcd)}\ _3\psi_3\left[\begin{array}{ccc} b/q,&c/q,&d/q\\[5pt]
q^2/b,&q^2/c,&q^2/d
\end{array}; q, \frac{q^4}{bcd} \right].
\end{align}
The above relation enables us to reduce the summation formula
\eqref{3psi3} to Jacobi's triple product identity.

\begin{exam}
A basic bilateral analogue of Dixon's sum \cite[Appendix
II.32]{Gasper-Rahman}:
\begin{align}
\label{4psi4}
&_4\psi_4\left[\begin{array}{cccc} -qa^{\frac{1}{2}},&b,&c,&d\\[5pt]
-a^{\frac{1}{2}},&aq/b,&aq/c,&aq/d
\end{array}; q, \frac{qa^{\frac{3}{2}}}{bcd} \right]\nonumber\\
&=\frac{(aq,aq/bc,aq/bd,aq/cd,qa^{\frac{1}{2}}/b,qa^{\frac{1}{2}}/c,qa^{\frac{1}{2}}/d,q,qa;q)_{\infty}}
{(aq/b,aq/c,aq/d,q/b,q/c,q/d,qa^{\frac{1}{2}},qa^{-\frac{1}{2}},qa^{\frac{3}{2}}/bcd;q)_{\infty}},
\end{align}
where $|qa^{\frac{3}{2}}/bcd|<1$.
\end{exam}

For the above formula, we may consider the substitution $d
\rightarrow d/q$ in (\ref{4psi4}) which suggests the iteration
relation
\begin{eqnarray}
\label{4psi4iterate}
& &_4\psi_4\left[\begin{array}{cccc} -qa^{\frac{1}{2}},&b,&c,&d\\[5pt]
-a^{\frac{1}{2}},&aq/b,&aq/c,&aq/d
\end{array}; q, \frac{qa^{\frac{3}{2}}}{bcd} \right] \nonumber \\[5pt]
&= &  \frac{(1-aq/bd)(1-aq/cd)(1-qa^{\frac{1}{2}}/d)}{(1-aq/d)(1-q/d)(1-qa^{\frac{3}{2}}/bcd)}\nonumber\\[5pt]
&  & \quad \times \, {}_4\psi_4\left[\begin{array}{cccc} -qa^{\frac{1}{2}},&b,&c,&d/q\\[5pt]
-a^{\frac{1}{2}},&aq/b,&aq/c,&aq^2/d
\end{array}; q, \frac{q^2a^{\frac{3}{2}}}{bcd} \right].
\end{eqnarray}
Let
\begin{equation}
F_k(a,b,c,d)=\frac{(-qa^{\frac{1}{2}},b,c,d;q)_k}{(-a^{\frac{1}{2}},aq/b,aq/c,aq/d)}\left(\frac{qa^{\frac{3}{2}}}{bcd}\right)^k
\end{equation}
and let
\begin{equation}
G_k(a,b,c,d)=\frac{(1-aq/bd)(1-aq/cd)(1-qa^{\frac{1}{2}}/d)}{(1-aq/d)(1-q/d)(1-qa^{\frac{3}{2}}/bcd)}F_k(a,b,c,d/q).
\end{equation}
By computation we obtain the Gosper pair:
\begin{eqnarray}\label{4psi4xk}g_k&=&\frac{-abdq^{k+1}-acdq^{k+1}+q^2a^\frac{3}{2}+a^2q^{k+2}-bcda^{\frac{1}{2}}q^{k+1}-da^{\frac{3}{2}}q^{k+1}+bcd^2q^{k}+bcd^2a^{\frac{1}{2}}q^{2k}}{(dq^k-q)(1+a^{\frac{1}{2}}q^k)(bcd-a^{\frac{3}{2}}q)} \nonumber\\[5pt]
&&\quad \times \, F_k(a,b,c,d),
 \\[5pt]
h_k&=&\frac{d(aq^k-c)(aq^k-b)}{(dq^{k-1}-1)(1+a^{\frac{1}{2}}q^k)(bcd-qa^{\frac{3}{2}})}F_k(a,b,c,d).\label{4psi4yk}
\end{eqnarray}
So  the iteration relation \eqref{4psi4iterate} holds.
From the symmetric property of the parameters $b,c$ and $d$, we have
\begin{align}
\label{4psi4iterate2}
&_4\psi_4\left[\begin{array}{cccc} -qa^{\frac{1}{2}},&b,&c,&d\\[5pt]
-a^{\frac{1}{2}},&aq/b,&aq/c,&aq/d
\end{array}; q, \frac{qa^{\frac{3}{2}}}{bcd} \right]\nonumber\\
&=\frac{(1-aq/bc)(1-aq^2/bc)(1-aq/bd)
(1-aq^2/bd)(1-aq/cd)(1-aq^2/cd)}{(1-aq/b)
(1-aq/c)(1-aq/d)(1-q/b)(1-q/c)(1-q/d)(1-qa^{\frac{3}{2}}/bcd)}\nonumber\\
&\qquad \times\frac{(1-qa^{\frac{1}{2}}/b)(1-qa^{\frac{1}{2}}/c)
(1-qa^{\frac{1}{2}}/d)}{(1-q^2a^{\frac{3}{2}}/bcd)(1-q^3a^{\frac{3}{2}}/bcd)}\nonumber
\\[5pt]
& \qquad \times {}_4\psi_4\left[\begin{array}{cccc} -qa^{\frac{1}{2}},&b/q,&c/q,&d/q\\[5pt]
-a^{\frac{1}{2}},&aq^2/b,&aq^2/c,&aq^2/d
\end{array}; q, \frac{q^4a^{\frac{3}{2}}}{bcd} \right].
\end{align}
By iteration, it follows that
\begin{align}
\label{4psi4last}
&_4\psi_4\left[\begin{array}{cccc} -qa^{\frac{1}{2}},&b,\ &c,&d\\[5pt]
-a^{\frac{1}{2}},&aq/b,&aq/c,&aq/d
\end{array}; q, \frac{qa^{\frac{3}{2}}}{bcd} \right]\nonumber\\
&=\frac{(aq/bc,aq/bd,aq/cd,qa^{\frac{1}{2}}/b,qa^{\frac{1}{2}}/c,qa^{\frac{1}{2}}/d;q)_{\infty}}
{(aq/b,aq/c,aq/d,q/b,q/c,q/d,qa^{\frac{3}{2}}/bcd;q)_{\infty}}H(a),
\end{align}
where
\begin{equation}
H(a)=\sum_{k=-\infty}^{\infty}\frac{(-qa^{\frac{1}{2}};q)_k}{(-a^{\frac{1}{2}};q)_k}q^{3{k
\choose 2}}\left(-qa^{\frac{3}{2}}\right)^k.
\end{equation}

Taking $b=-a^{\frac{1}{2}}$ and $c,d\rightarrow\infty$ in
\eqref{4psi4last} and by Jacobi's triple product identity
\eqref{jacobi}, it can be verified that
\begin{equation}
H(a)=\frac{(q,aq,q/a;q)_\infty}{(qa^{\frac{1}{2}},qa^{-\frac{1}{2}})_{\infty}},
\end{equation}
which leads to \eqref{4psi4}.

\begin{exam}
Bailey's very well-poised $_6\psi_6$-series identity
\eqref{6psi6}.
\end{exam}

Let us denote the $k$th term of (\ref{6psi6l}) as
\begin{equation}
\Omega_k(a;b,c,d,e)=\frac{(qa^{\frac{1}{2}}, -qa^{\frac{1}{2}}, b,
c, d, e;q)_k}{(a^{\frac{1}{2}}, -a^{\frac{1}{2}}, aq/b, aq/c,
aq/d, aq/e;q)_k}\left(\frac{qa^2}{bcde}\right)^k.
\end{equation}
Set $F_k=\Omega_k(a;b,c,d,e)$
and
\begin{equation}
G_k=\Omega_k(aq;b,c,d,eq)\times\frac{a(1-e)(1-aq)(bc-aq)(bd-aq)(cd-aq)}{(1-a)(b-aq)(c-aq)(d-aq)(bcde-a^2q)}.
\end{equation}
By computation, we find the following Gosper pair:
\begin{eqnarray}
g_k&=&\frac{a(bcdq^k-aq)(1-eq^k)}{(bcde-a^2q)(1-aq^{2k})}\Omega_k(a;b,c,d,e),
\nonumber \\[5pt]
h_k&=&\frac{(e-aq^k)(bcd-a^2q^{k+1})}{(bcde-a^2q)(aq^{2k}-1)}\Omega_k(a;b,c,d,e).
\end{eqnarray}
We note that the Abel pair derived from the above Gosper pair coincides with
the Abel pair given by Chu \cite{Chu}.

We next give the Gosper pair of a different
  iteration relation of the ${}_6\psi_6$
series by setting
$e\rightarrow e/q$ in (\ref{6psi6}), that is,
\begin{equation}\label{6psi6iterate}
\Omega(a;b,c,d,e)=\Omega(a;b,c,d,e/q)\times
\frac{(1-aq/be)(1-aq/ce)(1-aq/de)}{(1-aq/e)(1-q/e)(1-a^2q/bcde)}.
\end{equation}
We will see that the above iteration has the advantage that it
directly leads to the identity  \eqref{bcde} by taking into account
the symmetries in $b,c,d,e$. On the other hand, in this case the
Gosper pair does not have a simple expression.

Set $F_k=\Omega_k(a;b,c,d,e)$
and
\begin{equation}
G_k=\Omega_k(a;b,c,d,e/q)\times\frac{(1-aq/be)(1-aq/ce)(1-aq/de)}{(1-aq/e)(1-q/e)(1-a^2q/bcde)}.
\end{equation}

From the above iteration relation, we obtain the Gosper pair:
\begin{eqnarray}
g_k&=&\left(\frac{abceq^{k+1}+abdeq^{k+1}-a^2beq^{2k+1}+acdeq^{k+1}-a^2ceq^{2k+1}-a^2deq^{2k+1}-bcde^2q^{k}}{(bcde-qa^2)(eq^k-q)(aq^{2k}-1)}\right.\nonumber\\[5pt]
&&\left.+\frac{abcde^2q^{3k}-abcdeq^{2k+1}-a^2q^2+a^2eq^{k+1}+a^3q^{2k+2}}{(bcde-qa^2)(eq^k-q)(aq^{2k}-1)}\right)\Omega_k(a;b,c,d,e),\label{6psi6xk}\\[5pt]
h_k&=&\frac{qe(b-aq^k)(c-aq^k)(d-aq^k)}{(bcde-qa^2)(1-aq^{2k})(eq^k-q)}\Omega_k(a;b,c,d,e).\label{6psi6yk}
\end{eqnarray}
Since the parameters $b,c,d,e$ are symmetric in
\eqref{6psi6l}, we  obtain
\begin{align}
\Omega(a;b,c,d,e)&=\Omega(a;b/q,c/q,d/q,e/q)\nonumber\\
& \quad \times
\frac{(1-aq/bc)(1-aq^2/bc)(1-aq/bd)(1-aq^2/bd)}{(1-aq/b)(1-aq/c)(1-aq/d)(1-aq/e)}\nonumber\\
&\quad\times\frac{(1-aq/be)(1-aq^2/be)(1-aq/cd)(1-aq^2/cd)}
{(1-q/b)(1-q/c)(1-q/d)(1-q/e)}\nonumber\\
&\quad\times\frac{(1-aq/ce)(1-aq^2/ce)(1-aq/de)(1-aq^2/de)}
{(1-a^2q/bcde)(1-a^2q^2/bcde)(1-a^2q^3/bcde)(1-a^2q^4/bcde)}.\label{6psi62}
\end{align}

Again,  the limit
value can be given by Jacobi's triple product identity, so that we
arrive at \eqref{6psi6} in view of \eqref{6psi62}.

To conclude this section, we remark that our approach is feasible for the
computation of a Gosper pair for the $_{10}\psi_{10}$
 summation formula  of Chu \cite{Chu}.

\section{The Abel Pairs for Unilateral Summations}

The idea of Gosper pairs can be adapted to unilateral summation formulas
with a slight modification.
We also begin with an iteration relation guided by the closed product formula
which can be stated in the following form:
\begin{equation}
\sum_{k=0}^{\infty}F_k(a_1,a_2,\ldots,a_n)
=\sum_{k=0}^{\infty}G_k(a_1,a_2,\ldots,a_n).\label{abel5}
\end{equation}
We assume that $\sum_{k=0}^{\infty}F_k(a_1,a_2,\ldots,a_n)$
and $\sum_{k=0}^{\infty}G_k(a_1,a_2,\ldots,a_n)$ are convergent. Moreover, we
assume that the following limit condition holds:
 \begin{equation}
 \displaystyle{\lim_{k\rightarrow
\infty}}h_k=h_0.
\end{equation}

For the same reason as in the bilateral case, we see that
\[ F_k(a_1,a_2,\ldots,a_n)-G_k(a_1,a_2,\ldots,a_n)\]
 is a basic hypergeometric
terms so that we can resort to the $q$-Gosper algorithm to solve the
following equation:
\begin{equation}
F_k(a_1,a_2,\ldots,a_n)-G_k(a_1,a_2,\ldots,a_n)=h_{k+1}-h_k, \quad
k\geq 0.
\end{equation}
A Gosper pair $(g_k,h_k)$ is then given by
\begin{eqnarray*}
         g_k-h_k  & = &  F_k(a_1, a_2, \ldots, a_n), \\[5pt]
         g_k - h_{k+1} &   = &  G_k (a_1, a_2, \ldots, a_n).
\end{eqnarray*}
Therefore, one can use the Gosper pair $(g_k, h_k)$ to justify
(\ref{abel5}) by  the relation
\begin{equation}
\sum_{k=0}^{\infty}(g_k-h_k)= \sum_{k=0}^{\infty}(g_k-h_{k+1})
\end{equation}
and the limit condition ${\displaystyle\lim_{k\rightarrow \infty}}h_k=h_0$.

Given a Gosper pair it is easy to compute the corresponding
Abel pair which implies  iteration relation \eqref{abel5}
by the following unilateral Abel sum:
\begin{equation}
\sum_{k=0}^{\infty}A_k(B_k-B_{k-1})
=\sum_{k=0}^{\infty}B_k(A_k-A_{k+1}) ,\label{uAbellemma}
\end{equation}
which we call the unilateral Abel lemma.

The above approach is applicable to many classical unilateral
summation formulas including the $q$-Gauss sum, the $q$-Kummer
(Bailey-Daum) sum \cite[Appendix II.9]{Gasper-Rahman}, the $q$-Dixon
sum \cite[Appendix II.13]{Gasper-Rahman}, a $q$-analogue of Watson's
$_3F_2$ sum \cite[Appendix II.16]{Gasper-Rahman}, and a $q$-analogue
of Whipple's $_3F_2$ sum \cite[Appendix II.18]{Gasper-Rahman}, just
to name a few. Here we only give two examples to demonstrate this
technique.

\begin{exam}
The $q$-Gauss sum:
\begin{equation}
\label{2phi1}
_2\phi_1\left[\begin{array}{ccc} a,& b\\[5pt]
&c
\end{array}; q, \frac{c}{ab} \right] =  \frac{(c/a, c/b;q)_{\infty}}
{(c,c/ab;q)_{\infty}},
\end{equation}
\end{exam}
where $|c/ab|<1$.

Set
\begin{equation}
f(a,b,c)=\ _2\phi_1\left[\begin{array}{ccc} a,& b\\[5pt]
&c
\end{array}; q, \frac{c}{ab} \right].\label{2phi1left}
\end{equation}
Write the $k$th term of (\ref{2phi1left}) as
\begin{equation}
F_k(a,b,c)=\frac{(a, b;q)_k}{(q, c;q)_k}\left(\frac{c}{ab}
\right)^k.
\end{equation}
The iteration $c \rightarrow cq$ in (\ref{2phi1}) implies
\begin{equation}\label{2phi1iterate}
f(a,b,c)=\frac{(1-c/a)(1-c/b)}{(1-c)(1-c/ab)}f(a,b,cq).
\end{equation}
Let
\begin{equation}
G_k(a,b,c)=\frac{(1-c/a)(1-c/b)}{(1-c)(1-c/ab)}F_k(a,b,cq).
\end{equation}
Applying the $q$-Gosper algorithm to $F_k(a,b,c)-G_k(a,b,c)$, we
arrive at the Gosper pair:
\begin{eqnarray}
g_k & = & \frac{c-abq^k}{c-ab}F_k(a,b,c), \\[6pt]
h_k & = & \frac{ab(1-q^k)}{c-ab}F_k(a,b,c).
\end{eqnarray}
So we have the Abel pair:
\begin{eqnarray}
A_k & = &
\frac{(1-abq^k/c)}{(1-ab/c)}\frac{(a,b;q)_k}{(c,abq/c;q)_k}, \\[6pt]
B_k & = & \frac{(abq/c;q)_k}{(q;q)_k}\left(\frac{c}{ab}\right)^k.
\end{eqnarray}
Now we see that identity \eqref{2phi1} is true because of
 the unilateral Abel lemma \eqref{uAbellemma} and the
limit value $f(a,b,0)=1$.

\begin{exam}
The sum of Rogers' nonterminating very-well-poised $_6\phi_5$
series\cite[Appendix II.20]{Gasper-Rahman}:
\begin{eqnarray}
\lefteqn{
_6\phi_5\left[\begin{array}{cccccc} a,& qa^{\frac{1}{2}},
& -qa^{\frac{1}{2}},& b,& c,& d\\[5pt]
&a^{\frac{1}{2}},& -a^{\frac{1}{2}},& aq/b,& aq/c,& aq/d
\end{array}; q, \frac{aq}{bcd} \right] } \qquad \qquad \qquad
             \nonumber \\[6pt]
& = &   \frac{(aq, aq/bc, aq/bd, aq/cd;q)_{\infty}} {(aq/b, aq/c,
aq/d, aq/bcd;q)_{\infty}},   \label{6phi5}
\end{eqnarray}
where $|aq/bcd|<1$.
\end{exam}

Let us write
\begin{equation}
f(a,b,c,d)=\ _6\phi_5\left[\begin{array}{cccccc} a,& qa^{\frac{1}{2}},& -qa^{\frac{1}{2}},& b,& c,& d\\[5pt]
&a^{\frac{1}{2}},& -a^{\frac{1}{2}},& aq/b,& aq/c,& aq/d
\end{array}; q, \frac{aq}{bcd}\right].\label{6phi5left}
\end{equation}
Denote the $k$th term of (\ref{6phi5left}) by
\begin{equation}
F_k(a,b,c,d)=\frac{(a, qa^{\frac{1}{2}}, -qa^{\frac{1}{2}}, b, c,
d;q)_k}{(q, a^{\frac{1}{2}}, -a^{\frac{1}{2}}, aq/b, aq/c,
aq/d;q)_k}\left(\frac{aq}{bcd} \right)^k.
\end{equation}
The substitution $a \rightarrow aq$ in (\ref{6phi5}) leads to the iteration
relation
\begin{equation}\label{6phi5iterate}
f(a,b,c,d)=\frac{(1-aq)(1-aq/cd)(1-aq/bc)(1-aq/bd)}{(1-aq/b)(1-aq/c)(1-aq/d)(1-aq/bcd)}f(aq,b,c,d).
\end{equation}
Let
\begin{equation}
G_k(a,b,c,d)=\frac{(1-aq)(1-aq/cd)(1-aq/bc)(1-aq/bd)}{(1-aq/b)(1-aq/c)(1-aq/d)(1-aq/bcd)}F_k(aq,b,c,d).
\end{equation}
By computation we get the Gosper pair:
\begin{eqnarray}
g_k & = &
\frac{(1-aq^k)(q^k-aq/bcd)}{(1-aq^{2k})(1-aq/bcd)}F_k(a,b,c,d),\\[6pt]
h_k & = &
-\frac{(1-q^k)(1-a^2q^{k+1}/bcd)}{(1-aq^{2k})(1-aq/bcd)}F_k(a,b,c,d).
\end{eqnarray}
The corresponding Abel pair is as follows:
\begin{eqnarray}
A_k & = & \frac{(bcdq^k-aq)}{(bcd-aq)}
\frac{(b,c,d,a^2q^2/bcd;q)_k}{(aq/b,aq/c,aq/d,bcd/a;q)_k}, \\[6pt]
B_k & = &
\frac{(aq,bcd/a;q)_k}{(q,a^2q^2/bcd;q)_k}\left(\frac{aq}{bcd}\right)^k.
\end{eqnarray}
Therefore, the identity \eqref{6phi5} is a consequence of  the unilateral
Abel lemma \eqref{uAbellemma} and the limit value  $f(0,
b,c,d)=1$.

\vskip 5mm

\noindent
{\bf Acknowledgments.} This work was supported by the 973 Project on
Mathematical Mechanization, the Ministry of Education, the Ministry
of Science and Technology, and the National Science Foundation of China.

\end{document}